\documentclass[12pt]{article}

\usepackage[russian]{babel}
\usepackage{graphics}

\usepackage{epsfig}
\usepackage{cite}
\usepackage{enumerate}
\usepackage{amssymb}
\usepackage{amsmath}
\usepackage{graphicx}

\begin{document}


\noindent{\bf Solution to the Volterra operator equations of the 1st kind\\
with piecewise continuous kernels}\footnote{This work is supported by RFBR, project No. 11-08-00109,  Deutscher Akademischer Austauschdienst (DAAD), No.~A1200665. 
It is carried out within Federal Framework Programm ``Scientific and Scientific-Academic Staff of Innovative Russia'' within project No. 2012-1.2.2-12-000-1001-012.}\\

\noindent {\bf Denis Sidorov and Nikolai Sidorov}\\

\noindent{The sufficient conditions for existence and uniqueness of continuous
solutions of the
Volterra operator equations of the first kind with piecewise continuous kernel are derived.
The asymptotic
approximation of the parametric family of solutions are constructed
in case  of non-unique solution. The algorithm for
the solution's improvement is proposed using the successive approximations method.}

\noindent {\bf Keywords:} {\it Volterra operator equations of the first kind, asymptotic, discontinuous kernel, successive approximations, Fredholm's point, regularization.}\\

\noindent {\bf 1. Introduction}\\

Let us introduce area  $D = \{s,t; 0<s<t<T \}$ in the plane $s, t$ and
define continuous functions $s=\alpha_i(t), i=\overline{1,n},$ 
which have continuous derivatives for $t \in (0,T).$
Let us suppose that $\alpha_i(0)=0,$ $0<\alpha_1(t)< \dots <
\alpha_{n-1}(t) < t$  for $t \in (0,T),$
$0< \alpha_1^{\prime}(0)< \dots < \alpha_{n-1}^{\prime}(0)<1,$
and functions $s=\alpha_i(t), \, i=\overline{0,n},$
 $\alpha_0(t)=0, \, \alpha_n(t)=t,$
split the area $D$ into disjoint areas $D_1=\{s, t: 0\leq s <\alpha_1(t) \},$
$D_i = \{s,t : \alpha_{i-1}(t) < s < \alpha_i(t) , i=\overline{2,n}  \},$
$ \overline{D}={\bigcup\limits_1^n \overline{D_i}}.$
Let us introduce biparametric family of linear continuous operator-functions
  $K_i(t,s),$ defined for $t,s \in \overline{D}_i,$
$i=\overline{1,n},$
which are differentiable wrt $t$ and acting from  Banach space $E_1$ 
into Banach space $E_2$.
Therefore $K_i(t,s) \in \mathcal{L}(E_1 \rightarrow E_2)$
and $\frac{\partial K_i(t,s)}{\partial t} \in \mathcal{L}(E_1 \rightarrow E_2)$
for $t,s \in \overline{D_i}, i=\overline{1,n}.$
Let  the space of continuous functions $x(t)$ 
defined on  $[0,T]$ with ranges on  $E_1$ be denoted as $\mathbb{C}_{([0,T];E_1)}.$
Let us introduce the integral operator 
$$\int\limits_0^t K(t,s) u(s) ds \stackrel{\mathrm{def}}{=}  \sum\limits_{i=1}^n \int\limits_{\alpha_{i-1}(t)}^{\alpha_i(t)} K_i(t,s) x(s) ds 
\eqno{(1)}$$ 
with piecewise kernel
$$    K(t,s) = \left\{ \begin{array}{ll}
         \mbox{$K_1(t,s), \,\,\, t,s \in D_1$}, \\
         \mbox{\,\dots \,\,\,\,\, \dots \dots} \\
         \mbox{$K_n(t,s), \,\,\, t,s \in D_n$}. \\
        \end{array} \right. \eqno{(2)} $$
In this text we concentrate on the equation
$$\int\limits_0^t K(t,s)x(s)\,ds =f(t), \eqno{(3)} $$
where function $f(t)$ and function   $f^{\prime}(t)$ 
with ranges in $E_2$
 are defined and continuous for $t \in [0,T], \, f(0)=0.$
Let us call the equation (3) as the Volterra operator equation
of the first kind. The objective is solution to the equation (3)
in  $\mathbb{C}_{([0,T];E_1)}$.
Such problem and its numeric solution for $E_1=E_2=\mathbb{R}^1, n=1$
has been studied in number of papers, see e.g. [1--4].
Differentiation of the equation (3) lead to the integral-functional
equation and its solution is not unique in the general case. The reader 
may refer to  [5].  Therefore construction of the solution to equation (3)
can not be carried out with classic analytical mathods 
of the Volterra equations [6].
In this paper we address the equation (3) in the general case 
using the theory of operator equations with 
functionally perturbated argument of neutral type [7, \,8]. 
It is to be noted that the theory of the Volterra operator 
equations was initiated in  [9]. The principal results 
in this area with applications in integral geometry have been 
obtained in scientific school of M.M.Lavrentiev (see there references e.g. in 
the monograph [10]).

The paper is organized as follows. In Section 2
we continue our studies on the Volterra equations presented in papers [5,\, 11--17] 
and provide the sufficient conditions of existence and 
uniquness of continuous solution of the equation (3) with 
piece-wise continuous kernel (2).
The special case of such equations are the Volterra equations which model
the evolving dynamic systems [1--3]. To the best of our knowledge
such equations have not been studied in the literature. 
In the Section 2 the desired solution is constructed 
based on the ``step method" [8] from the theory of functional equations
combined with the method of successive approximations. 
In Section 3 and Section 4 we address the most interesting case when 
equation (3) has the family of solutions which depends on the free parameters.
 We propose the method for construction of asymptotic 
approximations of parametric solutions and we design the method
of the solution's refinement using the successive approximations.

In this work it is supposed that an operator  $K_n(t,t)$ 
has limited inverse operator for $t \in [0,T].$ Norms 
$||K_i(t,s)||_{\mathcal{L}(E_1\rightarrow E_2)},$ $ ||\frac{\partial K_i(t,s)}{\partial t}||_{\mathcal{L}(E_1\rightarrow E_2)},$  $i=\overline{1,n}$
of the linear operators are defined and they are continuous functions for 
 $t,s \in \overline{D}_i,$
$t\in [0,T].$ The conventional notations and terminology from [19, 20] are employed.\\

\noindent {\bf 2. Sufficient Conditions 
for Existence of the Unique Continuous Solution}\\

Since $f(0)=0$  differentiation of equation (3) 
leads to the  equivalent functional-operator equation
$$F(x)\stackrel{\mathrm{def}}{=} K_n(t,t)x(t) + \sum\limits_{i=1}^{n-1} \alpha_i^{\prime}(t)
\biggl\{ K_i(t,\alpha_i(t)) - $$
$$-K_{i+1}(t,\alpha_i(t)) \biggr\} x(\alpha_i(t)) + \sum\limits_{i=1}^n \int\limits_{\alpha_{i-1}(t)}^{\alpha_i(t)} K_i^{(1)}(t,s) x(s) \,ds - f^{\prime}(t)=0, \eqno{(4)} $$
where $\alpha_0=0, \, \alpha_n(t)=t.$

Let us introduce the function
$$D(t) \stackrel{\mathrm{def}}{=} \sum\limits_{i=1}^{n-1} | \alpha_i^{\prime}(t) | 
|| K_n^{-1}(t,t)
\{ K_i(t,\alpha_i(t)) - K_{i+1}(t,\alpha_i(t)) \} ||_{\mathcal{L}(E_1 \rightarrow E_1)}, $$
where $||\cdot||_{\mathcal{L}(E_1 \rightarrow E_1)}$ is the linear operator's norm.
Let the following condition be fulfilled:
\begin{enumerate} [{\bf (A)}]
\item $D(0) <1,$  $\sup\limits_{0<s<t<T} ||K_n^{-1}(t,t)K(t,s) ||_{\mathcal{L}(E_1 \rightarrow E_1)} \leq c <\infty.$
\end{enumerate}

Clearly the inequality $D(0)<1$ is fulfilled if $|\alpha_i^{(1)}(0)|$ are sufficiently small. 
Here and below an operator  $K(t,s)$  is defined with formula (2) in $\bigcup\limits_1^n D_i$. It's derivative wrt $t$ in common sense is defined as 
$$    K^{(1)}(t,s) = \left\{ \begin{array}{ll}
         \mbox{$K_1^{(1)}(t,s), \,\,\, t,s \in D_1$}, \\
         \mbox{\,\dots \,\,\,\,\, \dots \dots} \\
         \mbox{$K_n^{(1)}(t,s), \,\,\, t,s \in D_n$}. \\
        \end{array} \right.  \eqno{(5)} $$
for $t,s \in \bigcup\limits_1^n D_i.$

\noindent {\bf Theorem 1.} {\bf (Sufficient conditions for existence 
and uniqueness of the solution)}\\
{\it Let condition { (A)} be fulfilled, all the operators $K_i(t,s)$ in (2)
are continuous, and wrt $t$ they also have continuous derivatives, vector $f(t)$
has continuous derivatives, $f(0)=0.$
Then equation  (3) has unique solution in the class of continuous functions $\mathbb{C}_{([0,T];E_1)}$. And moreover the solution can be 
constructed using the step method combined with successive approximations.}\\

\noindent {\it Proof. }
Let us apply the operator $K_n^{-1}(t,t)$ to the equation (4). We get the following equation
$$
x(t) + A x + K x = \overline{f}(t),
\eqno{(6)}
$$
where the following notations are introduced\\
$A(t)x \stackrel{\mathrm{def}}{=} K_n^{-1}(t,t) \sum\limits_{i=1}^{n-1} \alpha_i^{(1)}(t) (K_i(t,\alpha_i(t)) - K_{i+1}(t,\alpha_i(t))) x(\alpha_i t),$\\
$Kx \stackrel{\mathrm{def}}{=} \sum\limits_{i=1}^n \int\limits_{\alpha_{i-1}(t)}^{\alpha_i (t)} K_n^{-1}(t,t) K_i^{(1)}(t,s) x(s)ds,\,\,$
$\overline{f}(t) \stackrel{\mathrm{def}}{=} K_n^{-1}(t,t) f^{(1)}(t).$
Let us fix $q<1$ and select  $h_1>0$ such as $ \max\limits_{0\leq t\leq h_1} D(t) = q<1.$
Such $h_1>0$ can be found because of condition { (A)}, and since the operator-functions $K_i(t,s)$ are continuous
in operator topology and because functions  $\alpha_i(t)$ are continue with their derivatives.
Lets $0<h<\min \{ h_1, \frac{1-q}{c} \},$
where constant value $c$ is defined in the condition {(A)}.
Let us split  $[0,T]$ into
$$
[0,h], \, [h, h+\varepsilon h], \, [h+\varepsilon h, h+ 2\varepsilon h], \dots .
\eqno{(7)}
$$
The contraction of the desired solution $x(t)$ onto the interval $[0,h]$ we denote as  $x_0(t),$  the contraction onto $$I_m = [(1+(m-1)\varepsilon)h, (1+m\varepsilon)h],\, m=1,2,\dots $$ we denote as $x_m(t).$ 
Let us select $\varepsilon$  from $(0,1]$ such as for $t\in I_m$ 
the ``perturbed'' arguments $\alpha_i(t) \in \bigcup\limits_{k=1}^{m-1} I_k, \, i=\overline{1,n-1}.$
If $0<\alpha_i^{(1)}(t) < \frac{1}{1+\varepsilon}$
for $t \in [0,T), \, i=\overline{1,n-1},$ then the mentioned above inclusion
is fulfilled on  $[0,T).$ 
Such inclusion enable the application of the ``step method'' (readers may refer to  [8],\, p. 199)) from the theory of functional-differential equations for construction of the solution  $x(t).$

For construction of the element  $x_0(t) \in \mathbb{C}_{([0,h],E_1)}$
we construct the sequence  $\{x_0^n(t)\}:$
$$x_0^n (t) = -A x_0^{n-1} - K x_0^{n-1} + \overline{f} (t), $$
$$x_0^0(t) = \overline{f}(t), \, t \in [0,h]. $$
Because of selection of $h$ we have an estimate  $||A+K||_{\mathcal{L}(\mathbb{C}_{([0,h];E_1)} \rightarrow \mathbb{C}_{([0,h];E_1)})} < 1.$

Therefore for $t\in [0,h]$ exists the unique solution $x_0(t)$ of the equation (6). The sequence 
$x_0^n(t)$ converge uniformly to this solution.
Let us continue the process of the desired solution 
construction for  $t\geq h,$ i.e. on 
$I_n, \, n=1,2, \dots .$
For sake of clarity let  $\varepsilon = 1$ in (7).

Then once we have the element $x_0(t) \in \mathbb{C}_{([0,h];E_1)}$ be constructed
we search for the element  $x_1(t) $ in the space $\mathbb{C}_{([h,2h];E_1)}$
of continuous vector-functions.
Let us find $x_1(t)$ from the Volterra equation of the 2nd kind
$$x(t) + \int\limits_h^t K_n^{-1}(t,t) K_t^{\prime}(t,s) x(s) \,ds =$$ $$= \overline{f}(t) - A x_0 - \int\limits_{0}^h K_n^{-1}(t,t)
K_t^{\prime}(t,s) {x}_0(s) \,ds$$
using the successive approximations. Here $x_0(h) = x_1(h).$

Let us introduce the continuous function
 $$
    \overline{x}_1(t) = \left\{ \begin{array}{ll}
         \mbox{$x_0(t), \,\,\,\, {0} \leq t \leq h$}, \\
         \mbox{$x_1(t), \,\, h \leq t \leq 2h $}, \\
\end{array}          
\right.
\eqno{(8)}
 $$
which is the contraction of the desired solution $x(t)$ onto $[0,2h].$
Then the element $x_2(t) \in$ $\mathbb{C}_{([2h,3h]; E_1)}$
can be computed with successive approximations from the Volterra
integral equation of the second kind
$$x(t) + \int\limits_{2h}^t K_n^{-1}(t,t) K_t^{\prime}(t,s) x(s) \, ds = $$ $$= \overline{f}(t) - A \overline{x}_1
- \int\limits_{0}^{2h} K_n^{-1}(t,t) K_t^{\prime}(t,s) \overline{x}_1(s)\, ds.$$
The desired solution $x(t) \in \mathbb{C}_{([0,T];E_1)}$ to the equation (3) can be constructed
on $N$th step, $N \geq \frac{T}{h}$.\\

\noindent {Example 1.}\\
Integral equation $$  \int\limits_{0}^{t/2} K_1(t-s) x(s) \,ds + \int\limits_{t/2}^{t} K_2(t-s) x(s)\, ds = f(t), \, 0 \leq t \leq T, $$ $f(0)=0,$
 $K_1(t-s) = K_2(t-s) +E, $ $K_1, K_2$ are matrices $m \times m,$ $E$ is the unit matrix,
$$|K_2^{-1}(0)|_{\mathcal{L}(\mathbb{R}^m \rightarrow \mathbb{R}^m)}<2,$$  matrix $K_2(t)$ and vector-function $f(t) = (f_1(t), \dots ,f_m(t))^{\prime}$
have continuous derivatives wrt $t,$ 
and satisfy conditions of the Theorem 1 and   has unique continuous
solution.\\

\noindent {Example 2.}\\
Boundary problem
$$\left\{ \begin{array}{ll}
         \mbox{$\int\limits_0^{t/2}  \biggl( \frac{\partial^2 x(t,y)}{\partial y^2} + x(t,y)  \biggr)\, dt + \int\limits_{t/2}^t \frac{\partial^2 x(t,y)}{\partial y^2} \,dt = f(t,y), \,\, 0\leq t \leq T,\,  0\leq y \leq 1$}, \\
         \mbox{$x(t,0) = 0, \,\,\,\,\,\, x(t,1)=0 $}, \\
\end{array}          
\right.$$
where function $f(t,y)$ is continuous wrt $y$ and has continuous derivative wrt $t,$  $f(0,y)=0,$
satisfy the conditions of the Theorem 1.
Desired continuous solution can be constructed as solution to the equivalent equation 
$$x(t,y) = -\frac{1}{2}\int\limits_0^t G(y,\xi) x(t,\xi) \,d\xi + \int\limits_0^1 G(y,\xi) 
f_t^{\prime}(t,\xi) \,d\xi $$
with contracting integral operator where
 $$
    G(y,\xi) = \left\{ \begin{array}{ll}
         \mbox{$(\xi-1)y, \,\,\,\, y \leq \xi $}, \\
         \mbox{$(y-1)\xi, \,\,\, \xi \leq y $} \\
\end{array}          
\right.
 $$
based on the successive approximations method.\\

\noindent {\bf 3. Construction of the asymptotic approximation $\hat{x}(t)$ 
of parametric family of solutions}\\

Let the following condition be fulfilled
\begin{enumerate} [{\bf (B)}]
\item Exist operator polynomials $\mathcal{P}_i = \sum\limits_{\nu+\mu = 1}^N K_{i \nu \mu} t^{\nu} s^{\mu}, \, i=\overline{1,n},$ where $K_{i,\nu,\mu} \in \mathcal{L}(E_1 \rightarrow E_2),$ are linear continuous operators,
vector-function $f^N(t) = \sum\limits_{\nu=1}^N f_{\nu} t^{\nu},$
polynomials $\alpha_i^N(t) = \sum\limits_{\nu = 1}^N \alpha_{i \nu} t^{\nu}, i=\overline{1,n-1},$
where
$0< \alpha_{11} < \alpha_{21} < \alpha_{21} < \dots < \alpha_{n-1,1} <1,$
are such as for $t\rightarrow +0,$ $s\rightarrow +0$ the following estimates
are fulfilled
$|| K_i(t,s) -\mathcal{P}_i(t,s) ||_{\mathcal{L}(E_1 \rightarrow E_2)} = \mathcal{O}((t+s)^{N+1}), \, i=\overline{1,n},$ 
$||f(t) - f^N(t)||_{E_2} = \mathcal{O}(t^{N+1}),$  $|\alpha_i(t) - \alpha_i^N(t)| = \mathcal{O}(t^{N+1}), i=\overline{1,n-1}$
 
\end{enumerate}
Expansion of powers  $t, s$ in the condition {(B)} we call as ``Taylor polynomials''
of the corresponding elements.

Let us introduce $j$-parametric family of linear operators
$$
B(j) = K_n(0,0) + \sum\limits_{i=1}^{n-1} (\alpha_i^{\prime}(0))^{1+j}(K_i(0,0)-K_{i+1}(0,0)),
$$
$j \in [0,\infty).$
Operator $B(j)$ which corresponds to the main ``functional'' part of the equation (4). 
We denote operator $B(j)$ as 
 {\it characteristic operator} of equation (4).

Let us follow the section 2 and consider equation (4) which is equivalent to the equation (3).
In contrast to the section 2 here we do not suppose that homogenius
equation for equation  (3) has only the trivial solution.
Therefore solution to the integral-functional
equation (4) can be non-unique. We follow  [5] and search
an asymptotic approximation of the particular solution
of non-homogenius equation (4)
as following polynomial
$$
\hat{x}(t) = \sum\limits_{j=0}^N x_j (\ln t) t^j.
\eqno{(9)}
$$

Let us demonstrate that coefficients $x_j$ with ranges in $E_1$  
in the general non-regular case depends on $\ln t$ and free parameters.
That is in line with that fact of existence of nontrivial
solution of the homogeneous equation.

There are regular and non-regular cases when it comes 
to determination of the coefficients $x_j$.\\

\noindent {\bf Definition 1.}
$j^*$ is the {\it regular point } of operator $B(j),$
 $B(j^*)$ has bounded inverse operator and {\it non-regular point} otherwise.\\

\noindent {\bf 3.1 Regular case: characteristic operator $B(j)$ has bounded 
inverse operator for $j\in (0,1, \dots ,N)$}\\

In this case the coefficients $x_j$ are constant vectors from 
$E_1$. Indeed, let us substitute (9) into the equation (4).
Then using the method of successive approximations
and condition {(B)} we get the recurrent sequence of linear
equations wrt $x_j:$
$$
B(0)x_0 = f^{\prime}(0),
\eqno{(10)}
$$
$$
B(j)x_j = M_j(x_0, \dots , x_{j-1}), \, j=1,\dots , N.
\eqno{(11)}
$$
Vector $M_j$ can be presented via solutions
$x_0, \dots , x_{j-1}$ and coefficients ``Taylor polynomials'' from the condition { (B)}.

Since operators $B(j)$ are invertible in the regular case then vectors
$x_0, \dots , x_N$ can be uniquely defined and therefore the asymptotic (9)
can be constructed.\\

\noindent {\bf 3.2 Non-regular case: operator $B(j)$ for   $j \in (0,1,\dots,N)$ 
has non-regular points}\\

\noindent Let us introduce definitions:

\noindent {\bf Definition 2.}
{\it Value $j^*$ is {\it simple singular Fredholm point of the operator} $B(j),$
if  $ B(j^*)$ is Fredholm operator ([20] p.219),
$\det \bigl[<B^{(1)}(j^*)\phi_i, \psi_k> \bigr]_{i,k=1}^r \neq 0,$
where $\{\phi_i \}_1^r$ is basis in $N(B(j^*)),$ $\{\psi_i\}_1^r$ is basis
in $N(B^{\prime}(j^*)),$ $B^{\prime}(j^*)$ is conjugate operator,
$B^{(1)}(j)$ is derivative wrt $j,$
computed with $j=j^*.$}\\

\noindent {\bf Definition 3.}
{\it Let $B(j^*)$ be Fredholm operator, $j^*$ we call singular Fredholm
point of index $k+1$ if $N(B(j^*)) \subset \bigcap\limits_{i=1}^k N(B^{(i)}(j^*)),$
$$\det \biggl[<B^{(k+1)}(j^*)\phi_i, \psi_k> \biggr]_{i,k=1}^r \neq 0,\, k\geq 1.$$}

It is to be noted that
$B^{(k)}(j) = \sum\limits_{i=1}^{n-1} (\alpha_i^{\prime}(0))^{1+j} a_i^k (K_i(0,0) - K_{i+1}(0,0)),$
where $a_i = \ln \alpha_i^{\prime}(0).$

\noindent {Remark 1.}
According to the definition 3 index of the simple Fredholm point is 1.
If $E_1=E_2=\mathbb{R}^1$ then $B(j)$ is function of $j.$ 
In this case the definition 2 means that 
$j^*$ is single root of the equation $B(j)=0,$ 
definition 3 means that  $j^*$ is $(k+1)$-multiple root of this equation.

Let us demonstrate that in non-regular case the coefficients $x_j$
will be polynomials of power  $\ln t$ and depends on 
arbitrary constants. Order of polynomials and number of 
arbitrary constants are connected with indices of singular
points of the operators
$B(j)$ and dimension $N(B(j))$.

Indeed since coefficient $x_0$ in the nonregular case 
can depends on  
$\ln t,$ then based on the method of undetermined
coefficients $x_0$ can be determined as the solution to 
the difference equation
  $$
K_n(0,0) x_0(z) + \sum\limits_{i=1}^{n-1} \alpha_i^{\prime}(0) (K_i(0,0) - K_{i+1}(0,0)) x_0(z+a_i) =
f^{\prime} (0),
\eqno{(12)}
$$
where $a_i = \ln \alpha^{\prime}(0), z = \ln t.$

Here there are three cases:

\noindent {\it  1st Case.}\\
Operator $B(0)$ has bounded inverse operator.
 Then coefficient $x_0$ does not depends on $z$
and can be determined uniquely from (10).

\noindent {\it 2nd Case.}\\
Let $j=0$ be the simple Fredholm point of the operator $B(j).$ Let us search the coefficient
 $x_0(z)$ from difference equation (12) as linear
vector-function 
$$
x_0(z) = x_{01} z + x_{02}. \eqno{(13)}
$$
We substitute (13) into (12) and we get the following two equations
for determination of vectors $x_{01}, x_{02}$ 
$$
B(0) x_{01} = 0,
\eqno{(14)}
$$
$$
B(0) x_{02} + B^{(1)}(0) x_{01} = f^{\prime}(0).
\eqno{(15)}
$$
Let $\{\phi_i\}_1^r$ be basis in $N(B(0)).$
Then $x_{01} = \sum\limits_{k=1}^r c_k \phi_k.$
Vector $c = (c_1, \dots , c_r)^{\prime}$
can be determined uniquely from the system
of linear algebraic equations
$$\sum\limits_{k=1}^r <B^{(1)}(0)\phi_k, \psi_i> c_k  = <f^{\prime}(0),\psi_i>, \, i=\overline{1,r}$$
with nonsingular matrix.
Further the coefficient $x_{02}$ can be determined from the equation (15)
with accuracy up to  $\mathrm{span}(\phi_1, \dots , \phi_r).$
For this reason we employ the formula
$$x_{02}=\sum\limits_{k=1}^r d_k \phi_k + \Gamma (f^{\prime}(0) - B^{(1)}x_{01}), $$
where $$\Gamma = (B(0)+ \sum\limits_{k=1}^r <\cdot,\gamma_k>z_k)^{-1} $$ is
Trenogin's regulariser (see [20], p.~221),
$d_1, \dots , d_r$ are arbitrary constants.
Therefore in the 2nd case the coefficient $x_0(z)$ 
is linear wrt
$z$ and depends on $r$ arbitrary constants.

\noindent {\it  3rd Case. } Let $j=0$ be the singular Fredholm point
of the operator $B(j)$ of index $k+1,$
 $k \geq 1.$ The solution $x_0(z)$ to the difference equation (12)
we search as polynomial
$$
x_0(z) = x_{01} z^{k+1} + x_{02} z^k + \dots + x_{0 k+1} z + x_{0 k+2}.
\eqno{(16)}
$$
Let us substitute of the polynomial (16)
into the system (12) and let us take into account the following identity $$\frac{d^k}{dj^k}B(j) = \sum\limits_{i=1}^{n-1} (\alpha_i^{\prime}(0))^{1+j} a_i^k (K_i(0,0) - K_{i+1}(0,0)),$$
where $a_i = \ln \alpha_i^{\prime}(0).$
Let us match the coefficients on powers of $z^{k+1}, z^k, \dots , z, z^0$
to zero. As result we get the sequence
of linear operator equations wrt  $x_{01}, x_{02}, \dots , x_{0 k+2}:$\\
$B(0) x_{01} = 0,$\\
$B(0) x_{02} + B^{(1)}(0) 
\biggl(\begin{array}{c}
	k+1\\ k
	\end{array}
	\biggr) x_{01} = 0,$\\
	$B(0) x_{0 l+1} + B^{(l)}(0) \biggl(\begin{array}{c}
	k+1\\ k+1-l
	\end{array}
	\biggr) x_{01} + B^{(l-1)}(0) \biggl(\begin{array}{c}
	k\\ k+1-l
	\end{array}
	\biggr) x_{02} + \dots  \\ \dots + B^{(1)}(0) 
	\biggl(\begin{array}{c}
	k+1-l+1\\ k+1-l
	\end{array}
	\biggr) x_{ol} = 0, \, l=1,\dots , k,
	$
$$
B(0) x_{0 k+2} + B^{(k+1)}(0) x_{01} + B^{(k)}(0) x_{02} + \dots B^{(1)}(0) x_{0 k+1} = f^{\prime}(0).
\eqno{(17)}
$$
Let us follow the definition 3 here. Therefore we have the following inclusion
$$
N(B(0)) \subset \bigcap_{i=1}^k N\biggl(\frac{d^i B(j)}{d j^i} \biggr|_{j=0}\biggr).
$$
Hence $B^{(i)}(0)x_{0 i+1}=0,\, i=\overline{0,k}$ and coefficients $x_{01}, \dots , x_{0 k+1}$
can be determined from  homogeneous equation $B(0)x=0$
 based on formulas $x_{0i} = \sum\limits_{j=1}^r c_{ij} \phi_j, \, i=\overline{1,k+1}.$
Therefore the equation (17) becomes
$$
B(0)x_{0,k+2} + B^{(k+1)}(0) x_{01} = f^{\prime}(0).
\eqno{(18)}
$$
Since $B(0)$ is the Fredholm operator and $\det \biggl[ <B^{(k+1)}(0)\phi_i, \psi_k> \biggr]_{i,k=\overline{1,r}} \neq 0,$
then vector $c^1 \stackrel{\mathrm{def}}{=} (c_{11}, \dots , c_{1r})^{\prime} $
can be determined uniquely from the conditions of 
resolvability of the equation (18).
Therefore we have
$$ x_{0\, k+2} = \sum\limits_{j=1}^r c_{{k+2}\, j} \phi_j + \hat{x}_{k+2}, $$
$\hat{x}_{k+2}$ is particular solution to the equation (18).
 Vector $c^{k+2} \stackrel{\mathrm{def}}{=}  (c_{k+2,1}, \dots , c_{k+2,r})^{\prime},$
as well as $c^i = (c_{i1}, \dots , c_{ir})^{\prime}, i=\overline{2,k+1},$
remains arbitrary. As result, the 3rd Case the coefficient $x_0(z)$
is polynomial of  $k+1$th power wrt $z$ and depends on 
$r \times (k+1) $ arbitrary constants.

Let us apply the method of undetermined 
coefficients and take into account the equality
$$\int t^j \ln^k t\, dt = t^{j+1} \sum\limits_{s=0}^k (-1)^s \frac{k(k-1) \dots (k-(s-1))}{(j+1)^{s+1}} \ln^{k-s} t.$$
We can construct the difference equations 
for determination of the coefficient
$x_1(z)$  ($z=\ln t$) and next coefficients 
of the asymptotic approximation (9).
Indeed let us take into account the definition
of the operator $F$ (see equation (4)). Then we get
$$
F(x)\biggr |_{x=x_0(z)+x_1(z) t}  = \biggl[ K_n(0,0)x_1(z) + \sum\limits_{i=1}^{n-1}
(\alpha_i^{\prime}(0))^2 (K_i(0,0) -
\eqno{(19)}
$$ 
$$-K_{i+1}(0,0)) x_1(z+a_i) + P_1(x_0(z)) \biggr] t + r(t), $$ with estimate $ r(t) = o(t).$
Here $P_1(x_0(z))$ is certain polynomial wrt $z.$ Its power is equal to
index of singular Fredholm point $j=0$ of the operator $B(j).$
From  (19) and because of the estimate $r(t)=o(t)$
for $t\rightarrow 0$ we conclude that the coefficient $x_1(z)$
have to satisfy the following difference equation
$$
K_n(0,0) x_1(z) + \sum\limits_{i=1}^{n-1} (\alpha^{\prime}(0))^2 \bigl(K_i(0,0) - K_{i+1}(0,0)\bigr) x_1(z+a_i)+
\eqno{(20)}
$$
$$+ P_1(x_0(z)) = 0. $$
If $j=1$ is regular point of the operator $B(j)$ then equation  (20)
has the solution $x_1(z)$ as polynomial of the same order as 
index of the singular Fredholm point $j=0$
of the operator $B(0).$
If $j=1$ also singular Fredholm point
of the operator $B(j)$ the the solution $x_1(z)$
can be constructed as polynomial of power  $k_0 +k_1,$
where $k_0$ and $k_1$ are indices of the corresponding  singular
Fredholm points $j=0,$ $j=1$ of the operator
$B(j)$. Coefficient $x_1(z)$ will depends on 
$r_0 k_0 + r_1 k_1$
arbitrary constants, where  $r_0 = \mathrm{dim} N(B(0)),$
$r_1 = \mathrm{dim} N(B(1)).$

Let us introduce the condition
\begin{enumerate} [{\bf (C)}]
\item
Let operator $B(j)$ for $j\in (0,1, \dots , N)$
has regular points only or singular fredholm
points  $j_1, \dots , j_{\nu}$
of indices $k_i,$ $dim N(Bj_i)=r_i, \, i=\overline{1,\nu}.$
\end{enumerate}

Then by the similar means we can determine 
the rest of the coefficients
$x_2(z), \dots , x_N(z)$ in $\hat{x}(t)$
from the sequence of difference equations
$$K_n(0,0)x_j(z) + \sum\limits_{i=1}^{n-1} (\alpha^{\prime}(0))^{1+j}\bigl(K_i(0,0)-K_{i+1}(0,0)\bigr)x_j(z+a_i)+$$
$$+\mathcal{P}_j(x_0(z), \dots , x_{j-1}(z)))=0, \, j=\overline{2,N}.$$
Therefore we have the following\\
\noindent {\bf Lemma 1.} \\
{\it Let conditions { (B)} and { (C)} be fulfilled. Then exists the vector-function
$\hat{x}(t) = \sum\limits_{i=0}^N x_i (\ln t) t^i,$ such as
$||F(\hat{x}(t))||_{E_2} = o(t^N),$ where operator $F$ is defined with formula (4).
The coefficients $x_i (\ln t)$ are polynomials of  $\ln t$
of increasing powers and they are smaller than $\sum\limits_j k_{j}$ 
singular Fredholm points $j\in\{0,1,2,\dots , N\}$ of characteristic operator $B(j)$.
Coefficients $x_i(\ln t)$ depends on $\sum_{j=0}^i \mathrm{dim} N(B(j)) k_j$
arbitrary constants.}\\

\noindent {\bf Remark 2.}
If $B(0)$ is the Fredholm operator and   $dim N(B(0)) \geq 1,$
then coefficient $x_0(\ln t)$ can be linear function of   $\ln t$ 
and vector-function $\hat{x}(t)$
will increase unboundly  for $t \rightarrow +0$
(briefly, $\hat{x} \in \mathbb{C}_{((0,T]; E_1)}$).\\

\newpage

\noindent {\bf 4. Theorem of existence of continuous
parametric solutions}\\

Since $0 \leq \alpha_i^{\prime}(0) <1, \, \alpha_i(0)=0, \, i=\overline{1,n-1},$
then for any $0<\varepsilon <1$
exists $T^{\prime} \in(0,T]$
such as $\max\limits_{i=\overline{1,n-1}, t\in [0,T^{\prime}]} |\alpha_i^{\prime} (t)| \leq \varepsilon$ and $\sup\limits_{i=\overline{1,n-1}, t\in (0,T^{\prime}]} \frac{\alpha_i(t)}{t} \leq \varepsilon.$

Let us introduce the condition

\begin{enumerate} [{\bf (D)}]
\item Let operator $K_n(t,t)$ has inverse bounded operator for  $t\in [0,T^{\prime}]$ 
and $N^*$ is selected to have the following inequality fulfilled
$$\sup\limits_{t\in (0,T^{\prime})} \varepsilon^{N^*} \sum\limits_{i=1}^{n-1} \bigl|\alpha_i^{(1)}(t)\bigr|  \biggl|\biggl|K_n^{-1}(t,t) \biggl(K_i(t,\alpha_i(t)) - K_{i+1}(t,\alpha_i(t))\biggr)\biggr|\biggl| _{\mathcal{L}(E_1 \rightarrow E_1)}\leq q <1.
$$
\end{enumerate}
 
\noindent {\bf Lemma 2.} \\
{\it Let condition { (D)} be fulfilled. Let in  $\mathbb{C}_{([0,T^{\prime}];E_1)}$
exists element $\hat{x}(t)$ such as for  $t\rightarrow +0$
$$||F(\hat{x}(t))||_{E_2} = o(t^N), N \geq N^*.$$ 
Then equation (3) in $\mathbb{C}_{([0,T^{\prime}];E_1)}$
has solution
$$
x(t) = \hat{x}(t) + t^{N^*} u(t),
\eqno{(22)}
$$
where $u(t)$ is  unique and defined with successive approximations.
}\\

\noindent  {\bf Proof. }
Let us substitute (22) into the equation  (4).
We get the following integral-functional equation for determination of  function
 $u(t)$
$$
K_n(t,t) u(t) + \sum\limits_{i=1}^{n-1} \alpha_i^{\prime}(t)
\biggl(\frac{\alpha_i(t)}{t} \biggr)^{N^*}  \biggl( K_i(t, \alpha_i(t))  -
\eqno{(23)}
$$
$$
-K_{i+1}(t,\alpha_i(t)) \biggr) u(\alpha_i(t)) +$$ $$+\sum\limits_{i=1}^n
\int\limits_{\alpha_{i-1}(t)}^{\alpha_{i}(t)} K_i^{(1)}(t,s) \biggl(\frac{s}{t}\biggr)^{N^*} u(s) \, ds
+ F(\hat{x}(t)) / t^{N^*} = 0.
$$
Let us introduce the linear operators
$$
Lu \stackrel{\mathrm{def}}{=}  K_n^{-1}(t,t)  \sum\limits_{i=1}^{n-1}
\alpha_i^{\prime}(t) \biggl(\frac{\alpha_i(t)}{t} \biggr)^{N^*}
\biggl \{ K_i(t, \alpha_i(t))  -$$
$$-K_{i+1}(t,\alpha_i(t)) \biggr\} u (\alpha_i(t)), $$
$$K u \stackrel{\mathrm{def}}{=} \sum\limits_{i=1}^n \int\limits_{\alpha_{i-1}(t)}^{\alpha_i(t)}
K_n^{-1}(t,t) K_i^{(1)}(t,s) (s/t)^{N^*} u(s)\, ds. $$
Then system (23) can be rewritten as follows
$$u+ (L+K) u  = \gamma(t), $$
where $\gamma(t) = K_n^{-1}(t,t) F(\hat{x}(t))/t^{N^*}$ 
is continuous vector-function.
Let us introduce the a Banach space $X$ of continuous functions $u(t)$
with ranges in a Banach space $E_1$
and following norm $$ ||u||_l = \max\limits_{0\leq t \leq T^{\prime}} e^{-lt} ||u(t)||_{E_1}, \, l>0.$$
Norm of linear functional operator $L$ satisfies the estimate
$$||L ||_{\mathcal{L}(X \rightarrow X)} \leq q < .$$
 because of inequality $\sup\limits_{t \in (0,T^{\prime}]} \frac{\alpha_i(t)}{t} \leq \varepsilon <1$
and based on the condition { (D)}  $\forall l \geq 0.$ 
Moreover for integral operator  $K$
the following estimate is fulfilled
$$||K ||_{\mathcal{L}(X \rightarrow X)} \leq q_1 < 1-q$$ for big enough $l.$
As results, for big enough  $l>0$ we have
$$||L+K ||_{\mathcal{L}(X \rightarrow X)} < 1,$$
i.e. the linear operator  $L+K$ is contracting in space $X.$
Therefore the sequence  $\{u_n \}$ converges. Here
 $u_n = -(L+K)u_{n-1} + \gamma(t), \, u_0=\gamma(t).$\\

   \noindent
{\bf Theorem 2. (Main Theorem)}\\
{\it Let conditions  { (B)}, { (C),} { (D)}, $f(0)=0$ be fulfilled. 
Let operator  $B(0)$ has bounded inverse operator.
Then equation (3) in space $\mathbb{C}_{([0,T];E_1)},$  
$0 \leq t\leq T^{\prime} \leq T$ has solution
$$ x(t) = \hat{x}(t) + t^{N^*} u(t), $$
which depends on $\sum\limits_{i=1}^{\nu} r_i k_i$
arbitrary constants.
Moreover, element $\hat{x}$ can be constructed as logarithm-power
sum (9), and $u(t)$
is uniquely computed with successive approximations.
The  asymptotic estimate  $||x(t) - \hat{x}(t)||_{E_1} = \mathcal{O}(t^{N^*})$
is fulfilled for $t\rightarrow +0.$}\\

   \noindent {\bf Proof. }
Let us employ Lemma 1 and take into account 
conditions of the theorem. Therefore the construction of 
asymptotic approximation
$\hat{x}(t)$ of the desired solution as following logarithm-power polynomial 
 $$\sum\limits_{i=0}^N x_i (\ln t) t^i$$ is possible.
And coefficients $x_i(\ln t)$ depend on certain number of 
arbitrary constants.
We can now apply Lemma 2 and  $x(t) = \hat{x}(t)+
t^{N^*} u(t).$ Therefore continuous function $u(t)$
can be constructed with successive approximations.
Theorem is proved.\\

The parametric family of solutions constructed on  $[0,T^{\prime}]$
we can continued with ``step method'' [8, p. 199] on the interval
$[0,T]$.\\
If  $j=0$ is Fredholm point of operator
 $B(j)$ and $ dim N(B(0)) \geq 1$ then based on Remark 3
the coefficient $x_0(\ln t)$ in the asymptotic $\hat{x}(t)$
can be linear function of $\ln t.$ In this case
 the solution $x(t) \in \mathbb{C}_{((0,T];E_1)}$
and grow unboundedly as $t\rightarrow +0$. \\

\noindent {\bf Example 3.}\\
Equation 
$$\int\limits_0^{t/2} \int\limits_0^{1}  K(y,y_1) x(s,y_1)\, dy_1 ds + \int\limits_{t/2}^t
\biggl( \int\limits_0^1 K(y,y_1) x(s,y_1)\, dy_1 - 2 x(s,y) \biggr)\, ds = g(y) t, $$
where $0<t<\infty, \, 0<y<1,$ 1 is eigen value of continuous symmetric
kernel $K(y,y_1)$
of rank $r,$ $\{\phi_1(y), \dots, \phi_r(y) \}$ is corresponding
orthonormal system of eigen functions on $[0,1]$, $g(y) \in \mathbb{C}_{[0,1]},$
meets the conditions of Theorem 2. In such case $j=0$ is simple
singular Fredholm point of corresponding characteristic operator
 $B(j).$ The solution of equation is following
 $$x(t,y) = -\frac{\ln t}{\ln 2} \sum\limits_{i=1}^r \int\limits_0^1 \phi_i(y) \phi_i(y_1) f(y_1)\, dy_1
+ c_1 \phi_1(y) + \dots +c_r \phi_r(y) + x_0(y), $$ $c_1, \dots , c_r$ 
are arbitrary constants,
$x_0(y)$ is the particular solution of Fredholm integral
equation of the second kind
$$x(y) = \int\limits_0^1 K(y,y_1)x(y_1)\, dy_1 - f(y) + \sum\limits_{i=1}^r \phi_i(y)
\int\limits_0^1 \phi_i(y_1) f(y_1)\, dy_1. $$


\noindent {\bf An Improvement of the Theorem 2. }

Let $\{j_1, \dots , j_{\nu}\} \in \mathbb{N} \bigcup\{0\}$  are Fredholm 
points of characteristic operator $B(j).$
 Let us construct generalized Jordan sets in sense of  [19, p.30]
for the characteristic operator 
$B(j)$. Then Theorem 2 can be improved.
Indeed let $j^*$ is singular Fredholm point of the operator $B(j).$
Let elements $\{\phi_i^{(l)} \}, i=\overline{1,r}, l=\overline{1,p_i}$ are constructed
and satisfy the following equalities:
$$B(j^*)\phi_i^{(1)} =0,$$
$$B(j^*)\phi_i^{(2)} + B^{(1)}(j^*) \bigl(\begin{smallmatrix} p_{i+1} \\ {p} \end{smallmatrix}\bigr) \phi_i^{(1)}= 0, $$
$$ \dots \dots \dots \eqno{(24)}$$
$$B(j^*)\phi_i^{(l+1)} + B^{(l)}(j^*) \bigl(\begin{smallmatrix} p_i+1 \\ {p_i+1-l} \end{smallmatrix}\bigr) \phi_i^{(1)} + \dots + B^{(1)}(j^*)\bigl(\begin{smallmatrix} p_i+1-l+1 \\ {p_i+1-l} \end{smallmatrix}\bigr)\phi_i^{(l)}= 0, $$
$$i=\overline{1,r}, \, l=\overline{1,p_i-1}. $$
Let also the following estimate be fulfilled
$$\det \biggl[ <B^{(p_i)}(j^*)\phi_i^{(1)} + B^{(p_i-1)}(j^*)\phi_i^{(2)}+ \dots + B^{(1)}(j^*)\phi_i^{(p_i)}, \psi_j> \biggr]_{i,j=\overline{1,r}} \neq 0,
\eqno{(25)}$$
where 
$\{ \phi_i \}_1^r$ is basis in $N(B^{\prime}(j^*)).$ 
Then we can follow known theory of Jordan sets of linear operators
( readers may here refer to [19, ch.30]) and  say 
that operator   $B(j)$ for $j=j^*$
has complete $B(j^*)$ --- Jordan set (CJS) $(\phi_i^{(l)} )_{i=\overline{1,r}, l=\overline{1,p_i}}.$
 We name $p_i$ as lengths of Jordan chains $$\stackrel{\longleftarrow \,\,\,\,\, p_i \,\,\,\,\, \longrightarrow }{\bigl(\phi_i^{(1)}, \dots , \phi_i^{(p_i)} \bigr)}, \, i=\overline{1,r}.$$

It is to be noted that CJS exists in $j^*$ if $j^*$ is singular Fredholm point
of index $p$  in sense of Definition 3. In this case $B^{(l)}(j^*)\phi_i^{(1)}=0, \,
i=\overline{1,r}, l=\overline{0,p-1}, $ Jordan chains
$$\stackrel{\longleftarrow \,\,\,\,\, p \,\,\,\,\, \longrightarrow }{\bigl(\psi_i^{(1)}, \dots , \psi_i^{(1)}\bigr)},\, i=\overline{1,r}$$ are stationary, i.e. they have the same lenght $p$, Condition (25) become following  $$\det [<B^{(p)}\phi_i^{(1)} , \psi_j>]_{i,j=\overline{1,r}} \neq 0. $$
Also condition { (C)} satisfied, i.e. we get the result of the Theorem 2.

Let us relax condition {(C)}:
\begin{enumerate} [{\bf (C1)}]
\item Characteristic operator $B(j)$ for $j \in (0,1, \dots ,N^*)$
has exactly $\nu$  singular Fredholm points $(j_1,\dots ,j_{\nu})$
with complete generalized Jordan sets, the rest of the values of this array are regular.
\end{enumerate}
It is to be noted that in condition (C1) Jordan chains
can be nonstationary.

If condition (C1) is satisfied then an asymptotic approximation $\hat{x}$ of desired
parametric family of solution of the equation (4) can be constructed.
Indeed let operator $B(0)$ is Fredholm operator, $\{\phi_i^{(1)} \}_{i=1}^r$ is the basis in
 $N(B(0)), \{ \phi_i^{(l)}\}_{i=\overline{1,r}, \, l=\overline{1,p_i}}$ 
is corresponding CJS satisfying equations (24) and condition (25)
for $j^*=0.$
Then the first coefficient $x_0(z)$ of  desired approximation $\hat{x},$
satisfy the difference equation (12) and can be constructed as polynomial
$$
x_0(z) = \sum\limits_{i=1}^r c_i \sum\limits_{l=1}^{p_i} \phi_i^{(l)} z^{p_i+1-l} +x_0,
\eqno{(26)}
$$  
where $c_1, \dots , c_r$ are constants and $x_0\in E_1$
to be determined. Substitution (26) in (12) yields linear equation wrt $x(0)$ 
$$ B(0)x_0 + \sum\limits_{i=1}^r c_i\bigl( B^{(p_i)}(0)\phi_i^{(1)} + B^{(p_i-1)}(0)\phi_i^{(2)} +
\dots  + B^{(1)}(0)\phi_i^{(p_i)} \bigr) = f^{\prime}(0). \eqno{(27)}$$
Because of inequality  (25) for $j^* = 0$ 
we can fine vector
 $(c_1, \dots , c_r)^{\prime}$ from the condition for
solvability of the equation (27). Next we can apply 
Trenogin regularizer ([20], \,p. 221)
 and construct the solution
$x_0$ of inhomogeneous equation (27) up to the basis $(\phi_1, \dots ,\phi_r)$
 in $N(B(0))$
Similarly, due to the condition {(C1)}
the rest of the coefficients  $x_2(z), \dots ,x_N(z)$
can be computed. Therefore based on the Lemma 2
we got Theorem 2 improved due to the condition (C) relaxation:\\

   \noindent  {\bf Theorem 3}.
Let conditions {(B)}, {(D)}, {(C1)}, $f(0)=0$ be satisfied.
Then equation (3) for  $0<t\leq T^{\prime}\leq T$
 has parametric family of solutions.\\

\noindent {\bf Conclusion}\\
In case of single equation when $E_1=E_2=\mathbb{R}^1,$ 
our method for solution of the difference systems solution
turn to be known A.O.Gelfond's method  (readers may refer to [18], \,p.338)
of construction of particular solutions of inhomogeneous
difference equations with polynomial right-hand side. 
We employed the results from functional analysis [15, \, 19,\, 20]
and the ideas of the Gelfond's method have been applied in the theory of 
the Volterra linear equations of the first kind with piece-wise continuous
operators acting on a Banach spaces. 
The constructed method 
can be employed for solution
of the class of integral-functional systems with singularities [14, \, 15].
If $f(0)\neq 0$ then the equation (1) does not have solution in class
of continuous functions. In this case
our method and results from the paper  [13] enable constrution
of the solutions in class of distributions [16,\, 17,\, 21].




\begin{thebibliography}{99}  

  \bibitem{lit1} Markova M. et al
  {About Models of Developing Systems of Glushkov Type and Their 
Applications in Electroenergetics,} 
Automationa and Remote Control, 7, 20--28 (2011).

  \bibitem{lit2} Yu. Yatsenko {Integral Models of Systems with Controllable Memory} (Naukova Dumka, Kiev, 1991).

  \bibitem{lit3} Apartsyn A.S. {Nonclassical Linear Volterra Equations of the First Kind } (De Gruyter, Walter, 2003).


  \bibitem{lit4}  Denisov A.M.,  Lorenzi A. {On a special Volterra integral equation of the first kind.} Boll. Un. Mat. Ital. B. Vol., {\bf (7)},  9, 443--457 (1995).

  \bibitem{lit5} Sidorov D.  {Volterra Equations of the First kind with Discontinuous Kernels in the Theory of Evolving Systems Control,} Studia Informatica Universalis, Vol.9, No.3, p.135--146 (2011).


    \bibitem{lit6} Marnitsky N.A. { Asymptotics of solutions to
the Volterra integral equations of the 1st kind}, DAN USSR, {\bf 269}, 1, 29--32 (1983).

 \bibitem{lit7} 
Sidorov N.A. and Trufanov A.V.  {Nonlinear operator equations with a functional perturbation of the argument of neutral type}, Differential Equations,  {\bf 45}, 12,  1840-1844 (2009).


\bibitem{lit8} 
Elsgoltz L.E. Qualitative Methods in Mathematical Analysis. Trans. Math. Mono, 12, American Math. Soc. 1964.

\bibitem{lit9} Lavrentiev M.M. and Buhgeim A.L. One class of operator equations of the first kind. Functional analysis and applications. 
Vol.7, No..4, 1973, p.44--53.

\bibitem{lit10} Lavrentiev M.M. and Saveliev L.Ya. Operator Theory and Ill-posed Problems.
IM SB RAS, Novosibirsk, 1999, 701p. (in Russian)


  \bibitem{lit11} Sidorov N.A. and Sidorov D.N. {Small solutions
of nonlinear differential equations in the neighborhood of branching points}, Izv VUZov. Mathematics, {\bf 5}, 53--61 (2011).

  \bibitem{lit12} Sidorov D.N., Sidorov D.N. {Convex majorants method in the theory of nonlinear Volterra equations,}   Banach J. Math. Anal.  {\bf 6},  1, 1 -- 10 (2012).





\bibitem{lit13} Sidorov N.A. and Sidorov D.N.  {Existence and construction of generalized solutions of nonlinear volterra integral equations of the first kind}, Differential Equations, {\bf 42},  9, 1312-1316 (2006).




Differential Equations, 2010, Volume 46, Number 6, Pages 882-891
  \bibitem{lit14} N. A. Sidorov, D. N. Sidorov and A. V. Krasnik {Solution of Volterra operator-integral equations in the nonregular case by the successive approximation method},
Differential Equations, 2010, Volume 46, Number 6, Pages 882-891.
 
 \bibitem{lit15} Sidorov N.A., Loginov B.V., Sinitsyn A.V., Falaleev  M.V. {Lyapunov-Schmidt methods in nonlinear analysis and applications}.
  Series on mathematics and its applications (Kluwer Academic Publishers, Dordrecht, 2002).


\bibitem{lit16} Sidorov D. {On impulsive control of nonlinear dynamical systems based on the Volterra series},  10th IEEE International Conference on  Environment and Electrical Engineering (EEEIC), 8-11 May 2011, Rome, Italy, 1--6 (2011). 

\bibitem{lit17} Sidorov D. and Sidorov N. {Generalized solutions
in problem of modeling of nonlinear dynamic systems with
the Volterra polynomials,} 
Automations and Remote Control, {\bf 6}, 127--132 (2011).







\bibitem{lit18} Gelfond A.O. {The calculus of finite differences}, Fizmatlit, Moscow, 1959, (in Russian).

\bibitem{li19} 
{M.M. Vainberg, V.A. Trenogin.} Theory of branching of solutions of non-linear equations, Noordhoff, 1974. (Translated from Russian).

\bibitem{lit20} { Trenogin V. A.} {Functional Analysis} (Fizmatlit, Moscow, 4th Ed., 2007).







\bibitem{21} Vladimirov V.S. Generalized functions in mathematical physics.
Nauka publ., Fizmatlit, Moscow, 1976 (in Russian).

\end{thebibliography}
\end{document}